\documentclass[12pt]{amsart}

\newcommand{\note}[2][\null]{% 
 \marginpar{ 
   \renewcommand{\baselinestretch}{.5} 
   \vspace{-1em}\hrule\vspace{3pt}% 
   \footnotesize\raggedright\textsf{#2\ifx#1\null\else\\\hfill--- 
     {\em #1}\fi}\vspace{1.5em} 
 }% 
} 

\usepackage{amsmath}
\usepackage{mathrsfs}
\usepackage{amsthm}
\usepackage{graphicx}
\usepackage{amsfonts}
\usepackage{wasysym}
\usepackage{textcomp}
\begin{document}
\title{An interesting application of Gegenbauer polynomials}
\author{Susanna Dann}

\address{Mathematics Department\\
Louisiana State University\\
Baton Rouge, Louisiana}
\email{sdann@math.lsu.edu}

\begin{abstract}
In this paper we will give a proof of $ \sum_{k=0}^m  \frac{\Gamma(\lambda+k) \Gamma(\lambda + m -k )}{\Gamma(\lambda) k! \Gamma(\lambda) (m-k)!} = \frac{\Gamma(m+2\lambda)}{\Gamma(2\lambda) m!} $ utilizing Gegenbauer polynomials.
\end{abstract}
\maketitle

\section{Preliminaries}
Gegenbauer polynomials belong to the family of orthogonal polynomials. As such they can be defined in different ways: as a solution of a certain differential equation, by a recursion relation or by means of a so called \text{\it{generating function}}. The last way is the most convinient for our purpose.

The coefficients $ C_m^{\lambda}(t) $ in the power series expansion of $ (1 - 2rt + r^2)^{-\lambda} $ for $ \lambda > 0 :$
\begin{equation}\label{eq1} 
(1 - 2rt + r^2)^{-\lambda} = \sum_{m=0}^\infty C_m^{\lambda}(t) r^m
\end{equation}
are called the \text{\it{Gegenbauer polynomials}} (\cite{leb72}, p.125).

Using $t=cos \varphi$ we can write $ (1 - 2rt + r^2)^{-\lambda} = [(1-r e^{i \varphi}) (1-r e^{- i \varphi})]^{-\lambda} $. Expanding the factors on the right-hand-side of the last equation we obtain:
\begin{equation}\label{eq2} 
\sum_{m=0}^\infty C_m^{\lambda}(t) r^m = \left[ \sum_{m=0}^\infty \binom{-\lambda}{m} e^{i m  \varphi} (-r)^m \right] \left[ \sum_{m=0}^\infty \binom{-\lambda}{m} e^{-i m  \varphi} (-r)^m \right].
\end{equation}

Note that $ \binom{-\lambda}{m} = (-1)^m \frac{\Gamma(\lambda + m)}{ \Gamma(\lambda) m!}$. We see that both series on the right-hand-side of (\ref{eq2}) converge absolutely for $ |r|<1 $ and uniformly in $\varphi$. Thus the left-hand-side of  (\ref{eq2}) converges absolutely for $ |r|<1 $ and uniformly in $t \in [-1, 1]$.

\section{Proof} 
We want to show that for any $ \lambda > 0 $
$$ \sum_{k=0}^m  \frac{\Gamma(\lambda+k) \Gamma(\lambda + m -k )}{\Gamma(\lambda) k! \Gamma(\lambda) (m-k)!} = \frac{\Gamma(m+2\lambda)}{\Gamma(2\lambda) m!}.$$

Let us evaluate (\ref{eq1}) for $t=1$.
\begin{align*}
     \sum_{m=0}^\infty C_m^{\lambda}(1) r^m &= (1 - 2r + r^2)^{-\lambda} \\
       &= (1 - r)^{-2\lambda} \\
       &= \sum_{m=0}^\infty \binom{-2\lambda}{m} (-r)^m.
  \end{align*}
Since $ \binom{-2\lambda}{m} = (-1)^m \frac{\Gamma(2\lambda + m)}{ \Gamma(2\lambda) m!}$, by coefficient comparison we obtain:
\begin{equation}\label{eq3} 
	C_m^{\lambda}(1) = \frac{\Gamma(2\lambda + m)}{ \Gamma(2\lambda) m!}.
\end{equation}

Recall that $ \left( \sum_{m=0}^\infty a_m r^m \right) \left( \sum_{m=0}^\infty b_m r^m \right) = \sum_{m=0}^\infty c_m r^m $ with $c_m =\sum_{k=0}^m a_k b_{m-k}$, if both series on the left-hand-side converge and at least one of them converges absolutely. We apply the last formula to $a_m = \frac{\Gamma(\lambda + m)}{ \Gamma(\lambda) m!}  e^{i m  \varphi}(t) $ and  $b_m = \frac{\Gamma(\lambda + m)}{ \Gamma(\lambda) m!}  e^{-i m  \varphi}(t) $ and obtain $$c_m = \sum_{k=0}^m \frac{\Gamma(\lambda + k) \Gamma(\lambda + m - k)}{ \Gamma(\lambda) k! \Gamma(\lambda) (m -k)!}  e^{i(2k - m) \varphi}. $$

Note that $c_m$ is nothing else but $ C_m^{\lambda}(t) = C_m^{\lambda}(cos \varphi) $. This gives
\begin{equation}\label{eq4} 
	C_m^{\lambda}(1) = \sum_{k=0}^m \frac{\Gamma(\lambda + k) \Gamma(\lambda + m - k)}{ \Gamma(\lambda) k! \Gamma(\lambda) (m -k)!}.
\end{equation}

Comparing (\ref{eq3}) and (\ref{eq4}) gives the result. \qed

\section*{Remark}
This proof was a by-prodct of my solution to an exersice, which asked to justify $ \frac{d}{dt} \sum_{m=0}^\infty C_m^{\lambda}(t) r^m =  \sum_{m=0}^\infty \frac{d}{dt} C_m^{\lambda}(t) r^m$.

\end{document}